\documentclass{amsart}
\usepackage{setspace,mathdefs,graphicx}

\newtheorem{definition*}{Definition}

\title{A new determinantal formula for the classical discriminant}
\author{Bradford Hovinen}
\address{
Warwick Mathematics Institute\\
Zeeman Building\\
University of Warwick\\
Coventry CV4 7AL United Kingdom}
\email{hovinen@math.utoronto.ca}

\begin{document}

\begin{abstract}
	According to several classical results by B\'ezout, Sylvester, Cayley, and others, the classical discriminant $D_n$ of degree $n$ polynomials may be expressed as the determinant of a matrix whose entries are much simpler polynomials in the coefficients of $f$. However, all of the determinantal formulae for $D_n$ appearing in the classical literature are equivalent in the sense that the cokernels of their associated matrices are isomorphic as modules over the associated polynomial ring. This begs the question of whether there exist formulae which are not equivalent to the classical formulae and not trivial in the sense of having the same cokernel as the $1\times 1$ matrix $(D_n)$.

	In this paper, we construct an explicit non-classical formula: the presentation matrix of the \textit{open swallowtail} first studied by Arnol'd and Givental. We study the properties of this formula, contrasting them with the properties of the classical formulae.
\end{abstract}

\maketitle

Let
$f(x,y):=a_0x^n+a_1x^{n-1}y+\cdots+a_{n-1}xy^{n-1}+a_ny^n$
be a homogeneous binary form of degree $n$ over an algebraically closed field $\KK$ of characteristic zero. Denote by $V$ the $\KK$-vector space of dimension $n+1$ with basis $\{a_0,\dots,a_n\}$. We may identify forms $f$, up to nonzero scalar multiple, with points in the projective space $\PProj(V)$. For $\alpha,\beta\in\KK$, not both zero, we say that the point $[\alpha:\beta]\in\PProj^1_{\KK}$ is a \textit{root of $f$ of multiplicity $k$} if the polynomial $(\beta x-\alpha y)^k$ divides $f$.

\begin{definition*}
	The \textit{(projective) classical discriminant of degree $n$ polynomials}, denoted $\Delta_n$, or $\Delta$ when $n$ is understood, is the locus of forms $f$ with a root of multiplicity at least two.
\end{definition*}

The variety $\Delta$ is a hypersurface in $\PProj(V)$ whose defining squarefree polynomial in $\KK[a_0,\dots,a_n]$ we call $D_n$. There is much interest in evaluating $D_n$ on a given binary form or on a family thereof. The direct approach --- constructing $D_n$ via, say, a parametrization of $\Delta$ and writing $D_n$ explicitly --- is infeasible for large values of $n$, since the number of terms of $D_n$ grows very quickly.\footnote{With respect to the $\ZZ^3$ grading $\deg a_i:=(1,n-i,i)$, $D_n$ is homogeneous of degree $(2(n-1),n(n-1),n(n-1))$. Thus one can compute an upper bound for the number of terms in $D_n$ which is exponential in $n$. There is no reason to believe that this bound is not tight, since it is unlikely that many of the coefficients of the terms of degree $(2(n-1),n(n-1),n(n-1))$ vanish.}

With this in mind, we are interested in the construction of a \textit{determinantal formula} for $D_n$, that is, a square matrix over $\KK[a_0,\dots,a_n]$ whose determinant is $D_n$. Of course, there exists a \textit{trivial} formula, namely, the $1\times 1$ matrix whose entry is $D_n$. More generally, we shall refer to any formula with matrix $A$ as \textit{trivial} if either $A$ or its classical adjoint is invertible.

There exist several nontrivial classical such formulae by such mathematicians as B\'ezout, Sylvester, and Cayley. However, they are all equivalent in the sense that they have locally isomorphic cokernels.\footnote{We omit a proof. See, e.g., \cite[Theorem 2.1.4]{Hov08} for a full account. Parts of the argument may be found in \cite[A.IV.78]{Bou81}.} Since one can easily construct arbitrarily many formulae by multiplying by invertible matrices on either side, and since the properties of such formulae are fundamentally the same, such a distinction is not very interesting for our purposes.

In this paper we construct a hitherto undiscovered determinantal formula for $D_n$ as the presentation matrix of a canonically defined $\Ocal_\Delta$-module called the \textit{open swallowtail}. This formula is nontrivial and inequivalent to the classical ones. As such, its associated matrix carries information which is not present in the classical formulae.

The remainder of this paper is organized as follows. Section \ref{SecPreliminaries} introduces some basic geometrical notation and facts about $\Delta$. Section \ref{SecAlgebraicDefinition} gives the main results of the paper and the definition of the open swallowtail. Section \ref{SecConstPresMatrix} describes in detail the construction of the presentation matrix of the open swallowtail. Section \ref{SecArnold} recalls Arnol'd's original geometric definition of the open swallowtail and shows that the algebraic and geometric definitions coincide. Section \ref{SecConductor} contains some technical results on the properties of the open swallowtail and Section \ref{SecApplication} contains a result on the data encoded by the matrix of the open swallowtail.

\section{\label{SecPreliminaries}Preliminaries}

\begin{figure}
	\begin{center}
		\includegraphics[width=6cm]{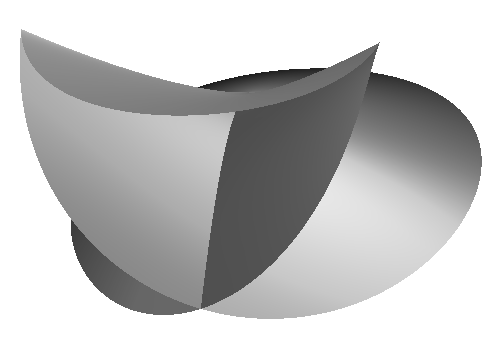}
	\end{center}
	\caption{\label{FigDiscriminant}The ``swallowtail''}
\end{figure}

Figure \ref{FigDiscriminant} shows the \textit{swallowtail}: a two-dimensional slice of the degree four discriminant $\Delta_4$ in the neighbourhood of the form $x^4$. The ``true'' discriminant near $x^4$ is the product of the surface shown in Figure \ref{FigDiscriminant} and an affine line. We see that $\Delta_4$ has a singular locus with two components, each in codimension one. The \textit{caustic}, denoted $\Gamma$, where the surface folds, is the locus of forms with a root of multiplicity three. The \textit{self-intersection locus}, where the surface crosses itself, is the locus of forms with more than one distinct pair of double roots. Both loci will be important in the sequel. These observations are true more generally, as the following theorem shows.

\begin{theorem}
	\label{ThmSingSigma}
	The singular locus of $\Delta$ is the locus of polynomials with either a root of multiplicity strictly greater than $2$ or more than one distinct pair of roots of multiplicity $2$. The former has codimension one and the latter has codimension one when $n\geq 4$ and is empty otherwise.

	Furthermore, the singular locus of $\Gamma$ is contained in its intersection with the self-intersection locus of $\Delta$.
\end{theorem}

\begin{proof}
	This is just a special case of \cite[Theorem 5.4]{Chi00}. The codimensions can be calculated by parametrizing the relevant loci to compute the dimensions thereof.
\end{proof}

In our study, the dictionary, introduced by Eisenbud in \cite{Eis80}, between maximal Cohen-Macaulay modules and matrix factorizations is crucial. Namely, given a determinantal formula $D_n=\det A$ for $D_n$, the pair $(A,\adj A)$, where $\adj A$ is the classical adjoint of $A$, is a matrix factorization of $D_n$ in the sense that $A(\adj A)=D_n\cdot I=(\adj A)A$, where $I$ is the identity matrix of appropriate size. The cokernel of $A$ is then a maximal Cohen-Macaulay module over $\KK[a_0,\dots,a_n]/(D_n)$. Furthermore, since $D_n$ is irreducible, given an MCM module $M$ over $\KK[a_0,\dots,a_n]/(D_n)$, a presentation of $M$ over $\KK[a_0,\dots,a_n]$ is a square matrix of determinant $D_n^r$, where $r$ is the rank of $M$. In particular, determinantal formulae correspond under this dictionary to \textit{maximal Cohen-Macaulay modules over $\KK[a_0,\dots,a_n]/(D_n)$ of rank one}. Under this dictionary, the classical determinantal formulae for $D_n$ correspond to the push-forward of the normalization $\bar{\Delta}$ of $\Delta$, which we now describe.

Denote by $W$ the $\KK$-vector space with basis $\{x,y\}$. Let $\Pbf:=\PProj(V)\times\PProj(W)$ and let 
\[F:=a_0x^n+a_1x^{n-1}y+\cdots+a_{n-1}xy^{n-1}+a_ny^n\in H^0(\Pbf,\Ocal(1,n))\]
be the universal homogeneous polynomial of degree $n$ in $x$ and $y$. Its partial derivatives\footnote{In general, for a polynomial $g(x,y)$, we denote $g_x:=\frac{\partial g}{\partial x}$ and $g_y:=\frac{\partial g}{\partial y}$.} $F_x:=\frac{\partial F}{\partial x}$ and $F_y:=\frac{\partial F}{\partial y}$ are sections of $\Ocal_{\Pbf}(1,n-1)$. Consider the \textit{incidence variety} $\bar{\Delta}$ defined by the sections $F_x$ and $F_y$. It is a smooth codimension two subvariety of $\Pbf$, as can be seen easily from the Jacobian criterion. The Euler identity \[ nF = xF_x+yF_y \] shows that on the affine pieces $U_y:=\{y\neq 0\}$ and $U_x:=\{x\neq 0\}$, respectively, $\bar{\Delta}$ coincides with the varieties defined by the sections $\{\frac{n}{y}F,F_x\}\subseteq\Gamma(U_y,\Ocal_{\Pbf}(1,n-1))$ and $\{\frac{n}{x}F,F_y\}\subseteq\Gamma(U_x,\Ocal_{\Pbf}(1,n-1))$. Thus points on $\bar{\Delta}$ are pairs $(f,t)\in\Pbf$ such that $t$ is a repeated root of $f$, the latter viewed as a homogeneous polynomial. In particular, the projection map $p_V:\Pbf\to\PProj(V)$, restricted to $\bar{\Delta}$, defines a map $\pi:\bar{\Delta}\to\Delta$.

\begin{proposition}
	\label{PropPHNormalization}
	The map $\pi:\bar{\Delta}\to\Delta$ described above is the normalization of $\Delta$.
\end{proposition}

\begin{proof}
	The map $\pi$ is finite since the number of preimages of a point on $\Delta$ is the number of distinct roots of multiplicity at least $2$ of the associated degree $n$ form, which is clearly a finite number. In addition, a generic polynomial of degree $n$ with a root of multiplicity $2$ has exactly one root of multiplicity exactly $2$, so $\pi$ is generically one-to-one. Finally, $\bar{\Delta}$ is smooth, hence normal. Thus the map $\pi$ is the normalization of $\Delta$, as claimed.
\end{proof}

We may make an analogous construction for the normalization of $\Gamma$. Namely, we define the incidence variety $\bar{\Gamma}$ to be the variety defined by the second-order partial derivatives $F_{xx}$, $F_{xy}$, and $F_{yy}$ of $F$, which are sections of $\Ocal_\Pbf(1,n-2)$. It is again easily seen that $\bar{\Gamma}$ is smooth. The same arguments as above show that the projection $p_V:\Pbf\to\PProj(V)$ maps $\bar{\Gamma}$ onto $\Gamma$ and that the restriction of $p_V$ to $\bar{\Gamma}$ is the normalization of $\Gamma$. The Euler identities $(n-1)F_x=xF_{xx}+yF_{xy}$ and $(n-1)F_y=xF_{xy}+yF_{yy}$ imply that $\bar{\Gamma}$ is a subvariety of $\bar{\Delta}$ of codimension one. We therefore have:

\begin{proposition}
	\label{PropGammaBarCM}
	The variety $\bar{\Gamma}$ embeds in $\bar{\Delta}$ as a smooth subvariety of codimension one. In particular, $\Ocal_{\bar{\Gamma}}$ is a Cohen-Macaulay module over $\Ocal_{\Delta}$ of codepth one.
\end{proposition}

\section{\label{SecAlgebraicDefinition}An algebraic definition}

In this section we develop an algebraic construction of the \textit{open swallowtail}, showing that it is a maximal Cohen-Macaulay module of rank one and therefore defines a determinantal formula for $D_n$.

The following theorem is key. It characterizes the sheaf of relative differentials $\Omega^1_{\bar{\Delta}/\Delta}$ of the normalization map $\bar{\Delta}\to\Delta$ described above and provides the main motivation for the algebraic definition of the open swallowtail, as well as the means to construct its presentation explicitly. Denote by $p$ the restriction of the projection $\Pbf\to\PProj(W)$ to $\bar{\Delta}$ and denote by $j$ the inclusion $\bar{\Gamma}\hookrightarrow\bar{\Delta}$.

\begin{theorem}
	\label{ThmCharOmega1}
	We have $\Omega^1_{\bar{\Delta}/\Delta}\cong j_*j^*p^*\Omega^1_{\PProj(W)/\KK}$ as $\Ocal_{\bar{\Delta}}$-modules.
\end{theorem}

\begin{proof}
	Let $\pi:\bar{\Delta}\to\Delta$ be the normalization map as given above. Let $i:\Delta\hookrightarrow\PProj(V)$ and $\bar{i}:\bar{\Delta}\hookrightarrow\Pbf$ be the natural embeddings. This gives rise to the diagram
	\begin{equation}
		\label{EqEmbeddingsDiagram}
		\begin{diagram}
			& & \bar{\Delta} & \rTo_\pi & \Delta \\
			& \ldTo_p & \dInto_{\bar{i}} & & \dInto_i \\
			\PProj(W) & \lTo_{p_W} & \Pbf & \rTo_{p_V} & \PProj(V).
		\end{diagram}
	\end{equation}

	We have $\Omega^1_{\Pbf/\KK}\cong p_V^*\Omega^1_{\PProj(V)/\KK}\oplus p_W^*\Omega^1_{\PProj(W)/\KK}$, where $p_W:\Pbf\to\PProj(W)$ and $p_V:\Pbf\to\PProj(V)$ are the natural projections. The Zariski-Jacobi sequence for the maps $\Pbf\to\PProj(V)\to\KK$ is just the split sequence associated to this direct sum decomposition:
	\begin{equation}
		\label{EqZJProduct}
		\begin{diagram}
			0 & \rTo & p_V^*\Omega^1_{\PProj(V)/\KK} & \rTo & \Omega^1_{\Pbf/\KK} & \rTo & \Omega^1_{\Pbf/\PProj(V)}\cong p_W^*\Omega^1_{\PProj(W)/\KK} & \rTo & 0.
		\end{diagram}
	\end{equation}

	The relations in \eqref{EqEmbeddingsDiagram}, along with the definition of $p$, imply that $\bar{i}^*p_V^*\Omega^1_{\PProj(V)/\KK}\cong\pi^*i^*\Omega^1_{\PProj(V)/\KK}$ and that $\bar{i}^*p_W^*\Omega^1_{\PProj(W)/\KK}\cong p^*\Omega^1_{\PProj(W)/\KK}$. Furthermore, because \eqref{EqZJProduct} is split, it remains exact after the application of $\bar{i}^*$. In view of the commutativity of \eqref{EqEmbeddingsDiagram}, we obtain the following commutative diagram:
	\[ \begin{diagram}
		& & \pi^*\Omega^1_{\Delta/\KK} & \rTo & \Omega^1_{\bar{\Delta}/\KK} & \rTo & \Omega^1_{\bar{\Delta}/\Delta} & \rTo & 0 \\
		& & \uTo & & \uTo \\
		0 & \rTo & \pi^*i^*\Omega^1_{\PProj(V)/\KK} & \rTo & \bar{i}^*\Omega^1_{\Pbf/\KK} & \rTo & p^*\Omega^1_{\PProj(W)/\KK} & \rTo & 0.
	\end{diagram} \]
	The top row is the Zariski-Jacobi sequence associated to the maps $\bar{\Delta}\to\Delta\to\KK$ and the bottom row is \eqref{EqZJProduct}. The vertical maps are the surjections $\pi^*i^*\Omega^1_{\PProj(V)/\KK}\to\pi^*\Omega^1_{\Delta/\KK}$ and $\bar{i}^*\Omega^1_{\Pbf/\KK}\to\Omega^1_{\bar{\Delta}/\KK}$ induced by the embeddings $i$ and $\bar{i}$. There is therefore an induced surjection $\rho:p^*\Omega^1_{\PProj(W)/\KK}\to\Omega^1_{\bar{\Delta}/\Delta}$, as in the following diagram:
	\begin{equation}
		\label{EqZJDiagram}
		\begin{diagram}[height=.5cm]
			& & 0 & & 0 & & 0 & & \\
			& & \uTo & & \uTo & & \uTo & & \\
			& & \pi^*\Omega^1_{\Delta/\KK} & \rTo & \Omega^1_{\bar{\Delta}/\KK} & \rTo & \Omega^1_{\bar{\Delta}/\Delta} & \rTo & 0 \\ \\
			& & \uTo & & \uTo & & \uDashto_{\exists\rho} & & \\ \\
			0 & \rTo & \pi^*i^*\Omega^1_{\PProj(V)/\KK} & \rTo & \bar{i}^*\Omega^1_{\Pbf/\KK} & \rTo & p^*\Omega^1_{\PProj(W)/\KK} & \rTo & 0.
		\end{diagram}
	\end{equation}

	It now suffices to prove that $\ker\rho=\Jcal\otimes p^*\Omega^1_{\PProj(W)/\KK}$, where $\Jcal$ is the ideal sheaf defining $\bar{\Gamma}$ in $\bar{\Delta}$. We have an injection $\Jcal\otimes p^*\Omega^1_{\PProj(W)/\KK}\hookrightarrow p^*\Omega^1_{\PProj(W)/\KK}$. Let $\Ical$ be the ideal sheaf defining $\bar{\Delta}$ in $\Pbf$. The cotangent sequence
	\[ \begin{diagram} \Ical/\Ical^2 & \rTo^d & \bar{i}^*\Omega^1_{\Pbf/\KK} & \rTo & \Omega^1_{\bar{\Delta}/\KK} & \rTo & 0 \end{diagram} \]
	then gives rise to the diagram
	\begin{equation}
		\label{EqCotDiagram}
		\begin{diagram}
			0 & & 0 \\
			\uTo & & \uTo \\
			\Omega^1_{\bar{\Delta}/\KK} & \rTo & \Omega^1_{\bar{\Delta}/\Delta} \\
			\uTo & & \uTo_\rho \\
			\bar{i}^*\Omega^1_{\Pbf/\KK} & \rTo & p^*\Omega^1_{\PProj(W)/\KK} \\
			\uTo & & \uTo \\
			\Ical/\Ical^2 & \rDashto & \ker\rho
		\end{diagram}
	\end{equation}
	whose middle rows come from \eqref{EqZJDiagram}. Thus there is an induced map $\Ical/\Ical^2\to\ker\rho$ making \eqref{EqCotDiagram} commute.

	For $i,j\geq 0$, let $F_{i,j}:=\frac{\partial^{i+j}}{\partial^i x\partial^j y}F$. Then $\Ical$ is defined by the sections $\{F_{1,0},F_{0,1}\}$ of $\Ocal(1,n-1)$. We claim that the composition $\Ical/\Ical^2\to\bar{i}^*\Omega^1_{\Pbf/\KK}\to p^*\Omega^1_{\PProj(W)/\KK}$ factors through a map $\chi:\Ical/\Ical^2\to\Jcal\otimes p^*\Omega^1_{\PProj(W)/\KK}$ given via
	\begin{align*}
		\chi(F_{0,1}) &:= F_{2,0}\ dx+F_{1,1}\ dy \\
		\chi(F_{1,0}) &:= F_{1,1}\ dx+F_{0,2}\ dy.
	\end{align*}

	The kernel of the surjection $\bar{i}^*\Omega^1_{\Pbf/\KK}\to\Omega^1_{\bar{\Delta}/\KK}$ is the subsheaf of $\bar{i}^*\Omega^1_{\Pbf/\KK}$ generated by the sections
	\begin{align*}
		\label{EqdF}
		dF_{1,0} &= \alpha_0\ da_0+\cdots+\alpha_n\ da_n + F_{2,0}\ dx+F_{1,1}\ dy \\
		dF_{0,1} &= \beta_0\ da_0+\cdots+\beta_n\ da_n + F_{1,1}\ dx+F_{0,2}\ dy,
	\end{align*}
	where $\alpha_0,\dots,\alpha_n,\beta_0,\dots,\beta_n$ are sections of $\Ocal(0,n-1)$. The summands $\alpha_0\ da_0+\cdots+\alpha_n\ da_n$ and $\beta_0\ da_0+\cdots+\beta_n\ da_n$ are in the component $\pi^*i^*\Omega^1_{\PProj(V)/\KK}$, while the summands $F_{2,0}\ dx+F_{1,1}\ dy$ and $F_{1,1}\ dx+F_{0,2}\ dy$ are in the component $p^*\Omega^1_{\PProj(W)/\KK}$. Since $\Jcal$ is generated by $F_{i,2-i}$ for $0\leq i\leq 2$, the image of the composition $\Ical/\Ical^2\to\bar{i}^*\Omega^1_{\Pbf/\KK}\to p^*\Omega^1_{\PProj(W)/\KK}$ lies in the image of the injection $\Jcal\otimes p^*\Omega^1_{\PProj(W)/\KK}\hookrightarrow p^*\Omega^1_{\PProj(W)/\KK}$, proving the claim.

	We now show that $\chi$ is surjective by a local calculation. Recall that, on $\bar{\Delta}$, for $i\in\{0,1\}$,
	\[ 0 = (n-1)F_{i,1-i} = xF_{i+1,1-i}+yF_{i,2-i}. \]
	On the affine piece $\{y\neq 0\}$, we therefore have
	\begin{align*}
		F_{i,2-i} &= -\frac{x}{y}F_{i+1,1-i} \\
		&= (-1)^{2-i}\left(\frac{x}{y}\right)^{2-i}F_{2,0}.
	\end{align*}
	In particular, since $\Omega^1_{\PProj(W)/\KK}$ is generated by $\frac{1}{y^2}(y\ dx-x\ dy)$ on the affine piece $\{y\neq 0\}$, $\Jcal\otimes p^*\Omega^1_{\PProj(W)/\KK}$ is generated by
	\[ F_{2,0}\frac{1}{y^2}(y\ dx-x\ dy) \]
	on this affine piece. On the other hand,
	\begin{eqnarray*}
		F_{i+1,1-i}\ dx+F_{i,2-i}\ dy
		&=& (-1)^{1-i}\left(\frac{x}{y}\right)^{1-i}F_{2,0}\ dx+(-1)^{2-i}\left(\frac{x}{y}\right)^{2-i}F_{2,0}\ dy \\
		&=& (-1)^{1-i}y\left(\frac{x}{y}\right)^{1-i}F_{2,0}\frac{1}{y^2}(y\ dx-x\ dy)
	\end{eqnarray*}
	Thus, taking $i=1$, the image of $\chi$ is also seen to be generated by $F_{2,0}\frac{1}{y^2}(y\ dx-x\ dy)$ on this affine piece, and $\chi$ is therefore surjective thereupon. The argument that $\chi$ is surjective on the other affine piece $\{x\neq 0\}$ is entirely symmetric.

	The arguments just given imply that the injection $\Jcal\otimes p^*\Omega^1_{\PProj(W)/\KK}\to p^*\Omega^1_{\PProj(W)/\KK}$ maps $\Jcal\otimes p^*\Omega^1_{\PProj(W)/\KK}$ into $\ker\rho:p^*\Omega^1_{\PProj(W)/\KK}\to\Omega^1_{\bar{\Delta}/\Delta}$. To prove the result, it remains to show that the image of $\Jcal\otimes p^*\Omega^1_{\PProj(W)/\KK}$ is in fact all of $\ker\rho$. To do this, we restrict to an affine open subset $U$ of $\bar{\Delta}$ on which $\Omega^1_{\bar{\Delta}/\Delta}$ is trivial. Being a quotient of the invertible module $p^*\Omega^1_{\PProj(W)}$, $\Omega^1_{\bar{\Delta}/\Delta}$ is a cyclic module on that affine piece and, applying the local trivialization identifies $\Gamma(U,\ker\rho)$ with an ideal $I\subseteq\Gamma(U,\Ocal_{\bar{\Delta}})$. The above arguments show that $I$ contains $\Gamma(U,\Jcal)$, which, since $\bar{\Gamma}$ is smooth and connected, is prime of codimension one. Hence any ideal properly containing $\Gamma(U,\Jcal)$ has codimension at least two.

	It follows from \cite[Proposition 5.1]{Chi00} that $\Omega^1_{\bar{\Delta}/\Delta}$ is at least supported at all points of $\Gamma$ outside of the self-intersection locus. Thus, as the set of such points is dense in $\Gamma$, it is supported at least on $\Gamma$, a codimension one set.
\end{proof}

Once we have characterized $\Omega^1_{\bar{\Delta}/\Delta}$, the following theorem, which motivates the definition of the open swallowtail, is immediate.

\begin{theorem}
	The universal derivation $d:\Ocal_{\bar{\Delta}}\to\Omega^1_{\bar{\Delta}/\Delta}$ is surjective.
\end{theorem}

\begin{proof}
	Locally, say, on the affine piece $\{a_0\neq 0, y\neq 0\}$, $\Omega^1_{\bar{\Delta}/\Delta}$ is cyclic and $d$ takes a local section $g\left(t,\frac{a_1}{a_0},\dots,\frac{a_n}{a_0}\right)\in\Ocal_{\bar{\Delta}}$ to $\frac{\partial g}{\partial t}\ dt$, where $t:=x/y$.
\end{proof}

\begin{definition}
	\label{DefOpenSwallowtail}
	The $n$-th \textit{(algebraic) open swallowtail} $\Scal_n$ is the kernel of $d:\Ocal_{\bar{\Delta}}\to\Omega^1_{\bar{\Delta}/\Delta}$. It is an $\Ocal_\Delta$-subalgebra of $\Ocal_{\bar{\Delta}}$. We refer $\Scal_n$ via $\Scal$ when the $n$ is understood.
\end{definition}

The following proposition shows that $\Scal$ indeed defines a determinantal formula for $D_n$.

\begin{proposition}
	\label{PropOSProps}
	The open swallowtail is a maximal Cohen-Macaulay module of rank one over $\Ocal_\Delta$.
\end{proposition}

\begin{proof}
	We have that $\Ocal_{\bar{\Delta}}$ is maximal Cohen-Macaulay over $\Ocal_\Delta$ and, by Proposition \ref{PropGammaBarCM} and Theorem \ref{ThmCharOmega1}, $\Omega^1_{\bar{\Delta}/\Delta}$ has codepth one. The exact sequence
	\[ \begin{diagram} 0 & \rTo & \Scal & \rTo & \pi_*\Ocal_{\bar{\Delta}} & \rTo & \pi_*\Omega^1_{\bar{\Delta}/\Delta} & \rTo & 0 \end{diagram} \]
	of $\Ocal_\Delta$ modules implies that the depth of $\Scal$ is at least the minimum of $\depth\Ocal_{\bar{\Delta}}$ and $\depth\Omega^1_{\bar{\Delta}/\Delta}+1$. Thus $\Scal$ is maximal Cohen-Macaulay. That $\Scal$ has rank one follows from its being embedded in $\Ocal_{\bar{\Delta}}$.
\end{proof}

We shall see that the formula defined by $\Scal$ is nontrivial and not equivalent to the classical formulae in Section \ref{SecConstPresMatrix}.

\section{\label{SecConstPresMatrix}Construction of the presentation matrix}

In this section we show how to construct a presentation matrix for $\Scal$ using the mapping cone construction applied to the short exact sequence
\[ \begin{diagram} 0 & \rTo & \Scal & \rTo & \pi_*\Ocal_{\bar{\Delta}} & \rTo^d & \pi_*\Omega^1_{\bar{\Delta}/\Delta} & \rTo & 0. \end{diagram} \]
We shall see in Section \ref{SecCayleyMethod} that $p_{V*}\Ocal_{\bar{\Delta}}$ and $p_{V*}\Omega^1_{\bar{\Delta}/\Delta}$, over $\PProj(V)$, have resolutions of the form
\[ \begin{diagram} 0 & \rTo & F_1 & \rTo^A & F_0 & \rTo & p_{V*}\Ocal_{\bar{\Delta}} & \rTo & 0 \end{diagram} \]
and, respectively,
\[ \begin{diagram} 0 & \rTo & G_2 & \rTo^{\partial_2} & G_1 & \rTo^{\partial_1} & F_0 & \rTo & p_{V*}\Omega^1_{\bar{\Delta}/\Delta} & \rTo & 0. \end{diagram} \]
We obtain the mapping cone diagram
\begin{equation} \label{EqMappingCone} \begin{diagram}
	0 & & 0 & & G_2 \\
	\dTo & & \dTo & & \dTo_{\partial_2} \\
	G_2 \oplus F_1 & \rTo & F_1 & \rTo_{D_1} & G_1 \\
	\dTo_B & & \dTo_A & & \dTo_{\partial_1} \\
	G_1 \oplus F_0 & \rTo & F_0 & \rTo_{D_0} & G_0 \\
	\dTo_E & & \dTo & & \dTo \\
	G_0 & & 0 & & 0.
\end{diagram} \end{equation}
The maps $D_0$ and $D_1$ are liftings of the universal derivation $d$ to the resolutions of $\Ocal_{\bar{\Delta}}$ and $\Omega^1_{\bar{\Delta}/\Delta}$. The map $B$ is given by the matrix
\[ \begin{gmatrix} \partial_2 & D_1 \\ 0 & A \end{gmatrix}, \]
while the map $E$ is given by the matrix
\[ \begin{gmatrix} \partial_1 & D_0 \end{gmatrix}. \]
Since $D_0$ is surjective, $E$ is surjective as well, so the complex
\[ \begin{diagram} 0 & \rTo & G_2\oplus F_1 & \rTo^B & G_1\oplus F_0 & \rTo^E & G_0 & \rTo & 0 \end{diagram} \]
is exact except at $G_1\oplus F_0$, where the homology is precisely $\Scal$. Therefore, we desire a presentation of the form
\[ \begin{diagram} G_2\oplus F_1 & \rTo^{\bar{B}} & G_1\oplus\ker D_0 & \rTo & \Scal & \rTo & 0, \end{diagram} \]
where $\bar{B}$ is given by the matrix
\[ \begin{gmatrix} \partial_2 & D_1 \\ 0 & \bar{A} \end{gmatrix}, \]
$\bar{A}$ being the restriction of $A$ to $A^{-1}(\ker D_0)$.

To construct the presentation matrix $\bar{B}$, we need three data:
\begin{itemize}
	\item the matrix $\bar{A}$,
	\item the matrix $\partial_2$, and
	\item the lifting $D_1$.
\end{itemize}
In Subsection \ref{SecCayleyMethod} we address the construction of $A$ and $\partial_2$, while in Subsection \ref{SecLiftingD} we address the construction of $D_1$. The restriction from $A$ to $\bar{A}$ is then quite easy and we address that at the end of this section.

\subsection{\label{SecCayleyMethod}Resolving $\Ocal_\Delta$ and $\Ocal_\Gamma$}

In this subsection, we construct resolutions of $\Ocal_{\bar{\Delta}}$ and $\Omega^1_{\bar{\Delta}/\Delta}$ --- or, rather, modules which are locally isomorphic to those. In view of Theorem \ref{ThmCharOmega1}, we replace the latter module with $\Ocal_{\bar{\Gamma}}$. Originally due to Cayley, the method we use was developed to construct determinantal formulae for the equation of the \textit{dual variety} $X^\vee$ of a given projective variety $X$. In our case, the variety $X$ is $\PProj(W)$, embedded via the $n$th Veronese embedding in $\PProj^n$, and $X^\vee$ is $\Delta_n$. Our treatment omits many of the technical points required for the general case. See \cite[Chapter 2]{GKZ94} for a thorough, modern treatment of this method.

\begin{figure}
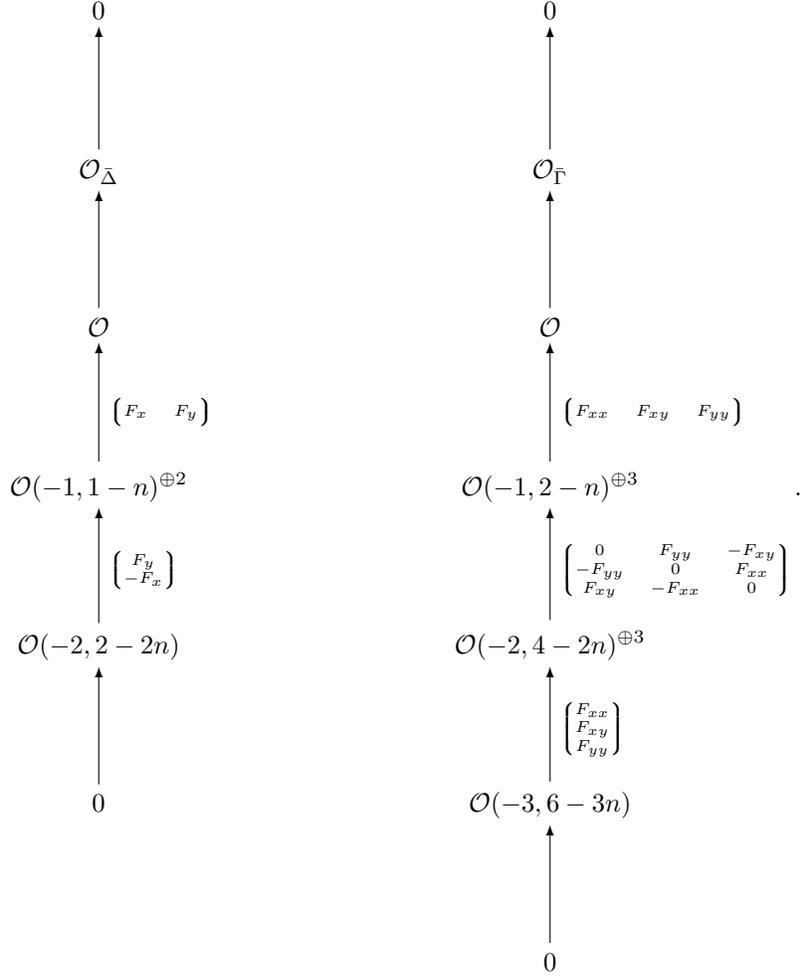

	\[ \begin{diagram}[width=6cm]
		0 & 0 \\
		\uTo & \uTo \\
		\Ocal_{\bar{\Delta}} & \Ocal_{\bar{\Gamma}} \\
		\uTo & \uTo \\
		\Ocal & \Ocal \\
		\uTo_{\mbox{\tiny $\begin{gmatrix} F_x & F_y \end{gmatrix}$}} & \uTo_{\mbox{\tiny $\begin{gmatrix} F_{xx} & F_{xy} & F_{yy} \end{gmatrix}$}} \\
		\Ocal(-1,1-n)^{\oplus 2} & \Ocal(-1,2-n)^{\oplus 3} \\
		\uTo_{\mbox{\tiny $\begin{gmatrix} F_y \\ -F_x \end{gmatrix}$}} & \uTo_{\mbox{\tiny $\begin{gmatrix} 0 & F_{yy} & -F_{xy} \\ -F_{yy} & 0 & F_{xx} \\ F_{xy} & -F_{xx} & 0 \end{gmatrix}$}} \\
		\Ocal(-2,2-2n) & \Ocal(-2,4-2n)^{\oplus 3} \\
		\uTo & \uTo_{\mbox{\tiny $\begin{gmatrix} F_{xx} \\ F_{xy} \\ F_{yy} \end{gmatrix}$}} \\
		0 & \Ocal(-3,6-3n) \\
		& \uTo \\
		& 0
	\end{diagram}. \]
	\caption{\label{FigKoszulComplexes} Koszul complexes resolving $\Ocal_{\bar{\Delta}}$ and $\Ocal_{\bar{\Gamma}}$}
\end{figure}

The normalizations $\bar{\Delta}$ and $\bar{\Gamma}$ of the discriminant $\Delta$ and respectively the caustic $\Gamma$ are global complete intersections in $\Pbf$, so their structure sheaves are resolved over $\Pbf$ via Koszul complexes, as shown in Figure \ref{FigKoszulComplexes}. To construct the required maps, we first twist the Koszul complex for $\bar{\Delta}$ by $\Ocal_\Pbf(0,n-2)$ and the Koszul complex for $\bar{\Gamma}$ by $\Ocal_\Pbf(0,n-3)$. We then construct the spectral sequence associated to the derived push-forward functor $Rp_{V*}$ and recover the required maps therefrom.

It follows from the projection formula that, for $l,p,q\in\ZZ$,
\[R^lp_{\PProj(V)*}\Ocal_{\Pbf}(p,q)\cong\Ocal_{\PProj(V)}(p)\otimes_\KK H^l(\PProj(W),\Ocal_{\PProj(W)}(q)).\]
In particular, $R^0p_{\PProj(V)*}\Ocal_{\Pbf}(p,q)=0$ when $q<0$ and $R^1p_{\PProj(V)*}\Ocal_{\Pbf}(i,j)$ when $q>-2$. The first page of the spectral sequence for $\Ocal_{\bar{\Delta}}$ is
\begin{align*}
	\Ocal_{\PProj(V)}\otimes_\KK H^0(\PProj(W),\Ocal_{\PProj(W)}&(n-2)) \\ \\
	& \Ocal_{\PProj(V)}(-2)\otimes_\KK H^1\left(\PProj(W),\Ocal_{\PProj(W)}(-n)\right),
\end{align*}
while the first page for $\Ocal_{\bar{\Gamma}}$ is
\begin{align*}
	\Ocal_{\PProj(V)}\otimes_\KK H^0(\PProj(W),\Ocal_{\PProj(W)}&(n-3)) \\ \\
	& \begin{diagram}
		\Ocal_{\PProj(V)}(-2)\otimes_\KK H^1\left(\PProj(W),\Ocal_{\PProj(W)}(1-n)^3\right) \\
		\uTo_{d_1^{-3,1}} \\
		\Ocal_{\PProj(V)}(-3)\otimes_\KK H^1\left(\PProj(W),\Ocal_{\PProj(W)}(3-2n)\right).
	\end{diagram}
\end{align*}

In both cases, the second page of the spectral sequence is
\[ \begin{diagram}
	\Ocal_{\PProj(V)}\otimes_\KK H^0\left(\PProj(W),\Ocal_{\PProj(W)}(n-k)\right) & & \\
	& \luTo_{d_2^{-2,1}} & \\
	& & \coker d_1^{-3,1},
\end{diagram} \]
where $k=2$ for $\bar{\Delta}$ and $k=3$ for $\bar{\Gamma}$. The cokernel of $d_2^{-2,1}$ is $\pi_*\Ocal_{\bar{\Delta}}(n-2)$, respectively $\pi_*\Ocal_{\bar{\Gamma}}(n-3)$, and the spectral sequence degenerates after this step. The map $A$ is then just a lifting of $d_2^{-2,1}$ to $\Ocal_{\PProj(V)}(-2)\otimes_\KK H^1\left(\PProj(W),\Ocal_{\PProj(W)}(-n)\right)$, while the map $\partial_2$ is just $d_1^{-3,1}$ in the diagram above.

We seek explicit formulae for $d_2^{-2,1}$ and $d_1^{-3,1}$. It is convenient to compute in the fibre over a fixed point $[a_0:\dots:a_n]\in\PProj(V)$. The restriction of the Koszul complex to this fibre is a map of vector spaces whose differentials vary polynomially in the coordinates $a_0,\dots,a_n$. In doing so, we replace the sheaves of the Koszul complex with \v{C}ech complexes which compute the cohomology of each term. We show the resulting double complex for both cases in Figure \ref{FigHCSpectralSequence}. Therein, $k=2$ for $\Ocal_{\bar{\Delta}}$ and $k=3$ for $\Ocal_{\bar{\Gamma}}$, and we use the convention that $\binom{k}{3}=0$ if $k=2$.

\begin{figure}
	\[ \begin{diagram}[height=1.4cm,loose]
		\begin{array}{c} \Gamma\left(U_x,\Ocal(n-k)\right) \\ \bigoplus \\ \Gamma\left(U_y,\Ocal(n-k)\right) \end{array} & \rTo_{d_H^{0,0}} & \Gamma\left(U_{xy},\Ocal(n-k)\right) \\
		\uTo_{d_V^{-1,0}} & & \uTo_{d_V^{-1,1}} \\
		\begin{array}{c} \Gamma\left(U_x,\Ocal(-1)^k\right) \\ \bigoplus \\ \Gamma\left(U_y,\Ocal(-1)^k\right) \end{array} & \rTo_{d_H^{-1,0}} & \Gamma\left(U_{xy},\Ocal(-1)^k\right) \\
		\uTo_{d_V^{-2,0}} & & \uTo_{d_V^{-2,1}} \\
		\begin{array}{c} \Gamma\left(U_x,\Ocal(k-n-2)^{\oplus\binom{k}{2}}\right) \\ \bigoplus \\ \Gamma\left(U_y,\Ocal(k-n-2)^{\oplus\binom{k}{2}}\right) \end{array} & \rTo_{d_H^{-2,0}} & \Gamma\left(U_{xy},\Ocal(k-n-2)^{\oplus\binom{k}{2}}\right) \\
		\uTo_{d_V^{-3,0}} & & \uTo_{d_V^{-3,1}} \\
		\begin{array}{c} \Gamma(U_x,\Ocal(2(k-n)-3)^{\oplus\binom{k}{3}}) \\ \bigoplus \\ \Gamma(U_y,\Ocal(2(k-n)-3)^{\oplus\binom{k}{3}}) \end{array} & \rTo_{d_H^{-k,0}} & \Gamma(U_{xy},\Ocal(2(k-n)-3)^{\oplus\binom{k}{3}}) \\
	\end{diagram} \]
	\caption{\label{FigHCSpectralSequence}Double complex for $\Kcal^k_\bullet$ over $[a_0:\cdots:a_n]$.}
\end{figure}

The differentials $d_1^{-i,0}$ and $d_1^{-i,1}$ are just the maps of the associated Koszul complex restricted to the given fibre. In particular, in the case of $\Ocal_{\bar{\Gamma}}$, the matrix of the last map $d_1^{-3,1}$ is a \textit{generalized Sylvester matrix} whose structure we now describe.

The map $d_1^{-3,1}$ is
\[ g\mapsto\left(\frac{\partial^2F}{\partial x^2}g,\frac{\partial^2F}{\partial x\partial y}g,\frac{\partial^2F}{\partial y^2}g\right), \]
where $g$ is a homogeneous element of $\KK[x^{-1},y^{-1}]$ of degree $3-2n<0$. The associated matrix is therefore divided horizontally into $3$ blocks of $n-2$ rows each. Each block is associated with a given second-order partial derivative of $F$. Each block is of the form
\[ \begin{gmatrix}
	\alpha_0 & \alpha_1 & \cdots & \alpha_{n-2} \\
	& \alpha_0 & \alpha_1 & \cdots & \alpha_{n-2} \\
	& & \ddots & \ddots & & \ddots \\
	& & & \alpha_0 & \alpha_1 & \cdots & \alpha_{n-2}
\end{gmatrix}, \]
where $\alpha_0,\dots,\alpha_{n-2}$ are the coefficients of the associated partial derivative of $F$. In particular, the entries are linear in the coefficients $a_0,\dots,a_n$.

We now characterize the map
\[A:H^1\left(\PProj(W),\Ocal(-n)\right)\to H^0\left(\PProj(W),\Ocal(n-2)\right).\]
It will be convenient to treat at the same time the corresponding map in the spectral sequence for $\bar{\Gamma}$, so we shall do so with the same notational conventions as above. Since $H^i(\PProj(W),\Ocal(-1))=0$ for all $i$, the map
\[d_H^{-1,0}:\Gamma\left(U_x,\Ocal(-1)^k\right)\oplus\Gamma\left(U_y,\Ocal(-1)^k\right)\to\Gamma\left(U_{xy},\Ocal(-1)^k\right)\] 
is an isomorphism. Thus the map
\[\left(d_H^{-1,0}\right)^{-1}:\Gamma\left(U_{xy},\Ocal(-1)^k\right)\to\Gamma\left(U_x,\Ocal(-1)^k\right)\oplus\Gamma\left(U_y,\Ocal(-1)^k\right)\]
is well-defined and $A$ is the composition
\begin{eqnarray*}
	\Gamma\left(U_{xy},\Ocal(k-n-2)^{\oplus\binom{k}{2}}\right) & \overset{d_V^{-2,1}}{\longrightarrow} & \Gamma\left(U_{xy},\Ocal(-1)^k\right) \\
	& \overset{\left(d_H^{-1,0}\right)^{-1}}{\longrightarrow} & \Gamma\left(U_x,\Ocal(-1)^k\right)\oplus\Gamma\left(U_y,\Ocal(-1)^k\right) \\
	& \overset{d_V^{-1,0}}{\longrightarrow} & \Gamma\left(U_x,\Ocal(n-k)\right)\oplus\Gamma\left(U_y,\Ocal(n-k)\right)
\end{eqnarray*}
restricted to $H^1\left(\PProj(W),\Ocal(k-n-2)^{\oplus\binom{k}{2}}\right)$.

In the case of $\Ocal_{\bar{\Delta}}$, we obtain the classical \textit{B\'ezout formula} for $D_n$, which we recall here. For a homogeneous binary form $f(x,y)$ of degree $n$, set
	\[ \Bcal_P(x_0,y_0,x_1,y_1):=\frac{f_x(x_0,y_0)f_y(x_1,y_1)-f_y(x_0,y_0)f_x(x_1,y_1)}{x_0y_1-x_1y_0}. \]
	Then $\Bcal_P(x_0,y_0,x_1,y_1)$ is a bihomogeneous polynomial which we write
	\[ \Bcal_P(x_0,y_0,x_1,y_1)=\sum_{i,j=0}^{n-2}b_{ij}x_0^iy_0^{n-i-2}x_1^jy_1^{n-j-2}. \]
	The \textit{B\'ezout matrix} is the matrix $B_P$ with entries $b_{ij}$. It gives a determinantal formula for $D_n$. The following proposition is proved in \cite{GKZ94}.

\begin{proposition}[\cite{GKZ94}, Chapter 2, Proposition 5.4]
	The matrix $B_P$ presents precisely the map $A$ as defined above.
\end{proposition}

\subsection{\label{SecLiftingD}Lifting the universal derivation}

We now describe how to construct the liftings $D_0$ and $D_1$ of the universal derivation $d$ to the resolutions of $\Ocal_{\bar{\Delta}}$ and $\Ocal_{\bar{\Gamma}}$. As in the previous, we fix a point $f\in\PProj(V)$ and work in the fibre of $p_V:\Pbf\to\PProj(V)$ over $f$.

There is one difficulty: while the universal derivation
\[ d:p_{V*}\Ocal_{\bar{\Delta}}\to p_{V*}\Omega^1_{\bar{\Delta}/\Delta} \]
is a well-defined map of $\Ocal_{\PProj(V)}$-modules, it cannot be twisted to form a map
\[ p_{V*}\left(\Ocal_{\bar{\Delta}}(n-2)\right)\to p_{V*}\left(\Omega^1_{\bar{\Delta}/\Delta}(n-3)\right). \]
We shall therefore construct maps $D_0$ and $D_1$ as in \eqref{EqMappingCone} and such that $D_0$ is a lifting of $d$ on the affine subset $U_y=\{y\neq 0\}$ of $\Pbf$. The mapping cone construction will then yield a module $\Scal'$ which is isomorphic to $\Scal$ on $U_y$, but may not agree with $\Scal$ globally. The affine piece $U_y$ corresponds to univariate polynomials of degree $n$, so this restriction suffices for our purposes.

As in the previous subsection, we shall fix a point $f\in\PProj(V)$ and compute in the fibre of $p_V$ over $f$. On $U_y$, $\Omega^1_{\PProj(W)/\KK}$ is trivial, being generated freely by $d(x/y)=\frac{1}{y^2}(y\ dx-x\ dy)$. In addition $\Ocal_{\PProj(W)}(i)$ is trivial for every $i$; the map $g\mapsto y^i g$ is an isomorphism $\Gamma\left(U_y,\Ocal_{\PProj(W)}\right)\to\Gamma\left(U_y,\Ocal_{\PProj(W)}(i)\right)$. We define
\[\tilde{D}_0:\Gamma\left(U_y,\Ocal_{\PProj(W)}(n-2)\right)\to\Gamma\left(U_y,\Omega^1_{\PProj(W)/\KK}(n-3)\right)\]
by composing with the local trivializations: for $h\in\Gamma(U_y,\Ocal_{\PProj(W)}(n-2))$, $\tilde{D}_0(h)=y^{n-3}d\left(y^{2-n}h\right)=h_x\ d(x/y)$. We then define
\[D_0:H^0\left(\PProj(W),\Ocal_{\PProj(W)}(n-2)\right)\to H^0\left(\PProj(W),\Omega^1_{\PProj(W)/\KK}(n-3)\right)\]
as the restriction of $\tilde{D}_0$ to $H^0\left(\PProj(W),\Ocal_{\PProj(W)}(n-2)\right)$.

Now we construct $D_1$. To do so, we lift $\tilde{D}_0$ through each map in the spectral sequence through which the homomorphisms
\[H^1\left(\PProj(W),\Ocal_{\PProj(W)}(1-n)^{\oplus 3}\right)\to H^0(\PProj(W),\Ocal_{\PProj(W)}(n-3))\]
and
\[H^1(\PProj(W),\Ocal_{\PProj(W)}(-n))\to H^0(\PProj(W),\Ocal_{\PProj(W)}(n-2))\]
are found.

We begin with the lifting
\[\begin{diagram}
	\Gamma\left(U_y,\Omega^1_{\PProj(W)/\KK}(-1)^{\oplus 3}\right) & \rTo & \Gamma\left(U_y,\Omega^1_{\PProj(W)/\KK}(n-3)\right) \\
	\uTo_{\tilde{D}_1^{(1)}} & & \uTo_{\tilde{D}_0} \\
	\Gamma\left(U_y,\Ocal_{\PProj(W)}(-1)^{\oplus 2}\right) & \rTo & \Gamma(U_y,\Ocal_{\PProj(W)}(n-2))
\end{diagram}.\]
Finding $\tilde{D}_1^{(1)}$ amounts to finding, for given $g_1,g_2\in\Gamma\left(U_y,\Ocal(-1)^{\oplus 2}\right)$, some sections $\tilde{g}_1,\tilde{g}_2,\tilde{g}_3\in\Gamma\left(U_y,\Ocal(-1)^{\oplus 3}\right)$ such that
\[(\tilde{g}_1f_{xx}+\tilde{g}_2f_{xy}+\tilde{g}_3f_{yy})\ d\left(\frac{x}{y}\right)=D_0(g_1f_x+g_2f_y).\]
Using the identities
\begin{align*}
	(n-1)f_x &= xf_{xx}+yf_{xy} \\
	(n-1)f_y &= xf_{xy}+yf_{yy},
\end{align*}
we compute directly
\begin{align*}
	(n-1)&\tilde{D}_0(g_1f_x+g_2f_y) \\
	=\ & (\left((n-1)g_1+xg_{1x}\right)f_{xx} + \left((n-1)g_2+yg_{1x}+xg_{2x}\right)f_{xy} \\
	& + yg_{2x}f_{yy})\ d\left(\frac{x}{y}\right).
\end{align*}
Thus
\begin{align*}
	(n-1)\tilde{g}_1 &= (n-1)g_1+xg_{1x}\\
	(n-1)\tilde{g}_2 &= (n-1)g_2+yg_{1x}+xg_{2x}\\
	(n-1)\tilde{g}_3 &= yg_{2x}.
\end{align*}
This provides the lifting $\tilde{D}_1^{(1)}$.

The next lifting
\[\begin{diagram}
	\Gamma\left(U_{xy},\Omega^1_{\PProj(W)/\KK}(-1)^{\oplus 3}\right) & \rTo & \Gamma\left(U_y,\Omega^1_{\PProj(W)/\KK}(-1)^{\oplus 3}\right) \\
	\uTo_{\tilde{D}_1^{(2)}} & & \uTo_{\tilde{D}_1^{(1)}} \\
	\Gamma\left(U_{xy},\Ocal_{\PProj(W)}(-1)^{\oplus 2}\right) & \rTo & \Gamma\left(U_y,\Ocal_{\PProj(W)}(-1)^{\oplus 2}\right)
\end{diagram}\]
is easy --- we just use the same definition as for $\tilde{D}_1^{(1)}$. We proceed to the final lifting
\[\begin{diagram}
	\Gamma\left(U_{xy},\Omega^1_{\PProj(W)/\KK}(1-n)^{\oplus 3}\right) & \rTo & \Gamma\left(U_{xy},\Omega^1_{\PProj(W)/\KK}(-1)^{\oplus 3}\right) \\
	\uTo_{\tilde{D}_1} & & \uTo_{\tilde{D}_1^{(2)}} \\
	\Gamma(U_{xy},\Ocal_{\PProj(W)}(-n)) & \rTo & \Gamma\left(U_{xy},\Ocal_{\PProj(W)}(-1)^{\oplus 2}\right)
\end{diagram}.\]
This is tantamount to finding, for a given section $g\in\Gamma(U_{xy},\Ocal_{\PProj(W)}(-n))$, sections \[h_1,h_2,h_3\in\Gamma\left(U_{xy},\Omega^1_{\PProj(W)/\KK}(1-n)^{\oplus 3}\right)\]
such that
\begin{align*}
	(n-1)(h_2f_{yy}-h_3f_{xy}) &= (n-1)g_1+xg_{1x} \\
	(n-1)(-h_1f_{yy}+h_3f_{xx}) &= (n-1)g_2+yg_{1x}+xg_{2x} \\
	(n-1)(h_1f_{xy}-h_2f_{xx}) &= yg_{2x},
\end{align*}
where $g_1:=f_yg$ and $g_2:=-f_xg$. We compute each in turn:
\begin{align*}
	(n-1)^2(h_2f_{yy}-h_3f_{xy})
	&= y\left((n-1)g+xg_x\right)f_{yy}+x\left(2(n-1)g+xg_x\right)f_{xy}, \\
	(n-1)^2(-h_1f_{yy}+h_3f_{xx})
	&= y^2g_xf_{yy}-x\left(2(n-1)g+xg_x\right)f_{xx}, \\
	(n-1)^2(h_1f_{xy}-h_2f_{xx})
	&= -y^2g_xf_{xy}-y\left((n-1)g+xg_x\right)f_{xx}.
\end{align*}
Thus we conclude
\begin{eqnarray}
	h_1 &=& -\frac{1}{(n-1)^2}y^2g_x \nonumber\\
	h_2 &=& \frac{1}{n-1}y\left(g+\frac{1}{n-1}xg_x\right) \label{EqLiftingD} \\
	h_3 &=& -\frac{1}{n-1}x\left(2g+\frac{1}{n-1}xg_x\right). \nonumber
\end{eqnarray}
Thus, for a section $g\in\Gamma(U_{xy},\Ocal_{\PProj(W)}(-n))$, we set
\[ \tilde{D}_1(g):=\left(-\frac{1}{(n-1)^2}y^2g_x,\frac{1}{n-1}y\left(g+\frac{1}{n-1}xg_x\right),-\frac{1}{n-1}x\left(2g+\frac{1}{n-1}xg_x\right)\right). \]
We then define $D_1$ as the restriction of $\tilde{D}_1$ to $H^1\left(\PProj(W),\Ocal_{\PProj(W)}(-n)\right)$.

The map $D_0$ is essentially differentiation with respect to $x$. The kernel of $D_0$ is therefore the sub-vector space of $F_0$ generated by $y^{n-2}$.

\begin{theorem}
	The matrix
	\begin{equation}
		\label{EqPresentOS}
		\begin{gmatrix} \partial_2 & D_1 \\ 0 & A_{(y^{n-2})} \end{gmatrix},
	\end{equation}
	where $\partial_2$ and $D_1$ are as given above and $A_{(y^{n-2})}$ is the row of the B\'ezout matrix $A$ corresponding to the basis element $y^{n-2}$ of $H^0(X,\Ocal(n-2))$, presents a module isomorphic on the open affine set $U_y$ to the open swallowtail $\Scal_n$.
\end{theorem}

We may also use this construction to count the minimal number of generators of $\Scal_n$ on the affine piece $\{a_0\neq 0\}$

\begin{proposition}
	\label{PropMinGenerators}
	On the affine piece $\{a_0\neq 0\}$, the open swallowtail $\Scal_n$ is minimally generated by $n-2$ elements. In particular, it is not isomorphic to the normalization of $\Delta_n$ and, when $n>3$, not free.
\end{proposition}

\begin{proof}
	After restricting to the affine piece $\{a_0\neq 0\}$, the now invertible variable $a_0$ appears in exactly $n-2$ of the $2(n-2)+(n-1)=3n-5$ columns of \eqref{EqPresentOS}. Furthermore, there are invertible elements in each of the $n-1$ columns of $D_1$. None of the other columns has an invertible element. Thus a minimal presentation matrix of $\Scal_n$ is of size $(n-2)\times (n-2)$. On the other hand, since the B\'ezout formula contains no units and is $(n-1)\times (n-1)$, the normalization $\bar{\Delta}$ is minimally $n-1$ generated.
\end{proof}

\begin{example}
	We construct a presentation matrix of the open swallowtail for the degree 4 discriminant.  In accordance with the above discussion, we construct a presentation of the open swallowtail over the affine piece $\{y\neq 0\}$. The matrix representing $\partial_2$ is
	\[\partial_2=\begin{gmatrix}
		12a_0 & 6a_1 & 2a_2 & 0 \\
		0 & 12a_0 & 6a_1 & 2a_2 \\
		3a_1 & 4a_2 & 3a_3 & 0 \\
		0 & 3a_1 & 4a_2 & 3a_3 \\
		2a_2 & 6a_3 & 12a_4 & 0 \\
		0 & 2a_2 & 6a_3 & 12a_4
	\end{gmatrix}.\]
	The map $D_1$ is given by the formulae \eqref{EqLiftingD} above. We compute
	\begin{eqnarray*}
		D_1(x^{-3}y^{-1}) &=& \left(-\frac{1}{3}yx^{-4},0,-\frac{1}{3}x^{-2}y^{-1}\right) \\
		D_1(x^{-2}y^{-2}) &=& \left(\frac{2}{9}x^{-3},\frac{1}{9}x^{-2}y^{-1},-\frac{4}{9}x^{-1}y^{-2}\right) \\
		D_1(x^{-1}y^{-3}) &=& \left(\frac{1}{9}x^{-2}y^{-1},\frac{2}{9}x^{-1}y^{-2},-\frac{5}{9}y^{-3}\right).
	\end{eqnarray*}
	Thus $D_1$ has the form
	\[D_1=\begin{gmatrix}
		0 & 0 & \frac{1}{9} \\
		0 & 0 & 0 \\
		0 & \frac{1}{9} & 0 \\
		0 & 0 & \frac{2}{9} \\
		-\frac{1}{3} & 0 & 0 \\
		0 & -\frac{4}{9} & 0
	\end{gmatrix}.\]

	The matrix $A$ is the B\'ezout matrix
	\[A=\begin{gmatrix}
		a_1a_3-16a_0a_4 & 2a_2a_3-12a_1a_4 & 3a_3^2-8a_2a_4 \\
		2a_1a_2-12a_0a_3 & 4a_2^2-8a_1a_3-16a_0a_4 & 2a_2a_3-12a_1a_4 \\
		3a_1^2-8a_0a_2 & 2a_1a_2-12a_0a_3 & a_1a_3-16a_0a_4
	\end{gmatrix}.\]
	Putting these together, $\Scal'$ is presented by
	\[\Scal'=\coker\begin{gmatrix}
		12a_0 & 6a_1 & 2a_2 & 0 & 0 & 0 & \frac{1}{9} \\
		0 & 12a_0 & 6a_1 & 2a_2 & 0 & 0 & 0 \\
		3a_1 & 4a_2 & 3a_3 & 0 & 0 & \frac{1}{9} & 0 \\
		0 & 3a_1 & 4a_2 & 3a_3 & 0 & 0 & \frac{2}{9} \\
		2a_2 & 6a_3 & 12a_4 & 0 & -\frac{1}{3} & 0 & 0 \\
		0 & 2a_2 & 6a_3 & 12a_4 & 0 & -\frac{4}{9} & 0 \\
		0 & 0 & 0 & 0 & \eta_1 & \eta_2 & \eta_3
	\end{gmatrix},\]
	where
	\begin{align*}
		\eta_1 &:= a_1a_3-16a_0a_4, \\
		\eta_2 &:= 2a_2a_3-12a_1a_4, \\
		\eta_3 &:= 3a_3^2-8a_2a_4.
	\end{align*}
	After applying necessary row operations, we obtain the following minimal presentation matrix:
	\[\Scal'=\coker\bgroup\begin{gmatrix}
		0 & 12a_0 & 6a_1 & 2a_2 \\
		-24a_0 & -9a_1 & 0 & 3a_3 \\
		12a_1 & 18a_2 & 18a_3 & 12a_4 \\
		\gamma_1 &
		\gamma_2 &
		\gamma_3 &
		0 \\
	\end{gmatrix}\egroup,\]
	where
	\begin{eqnarray*}
		\gamma_1 &=& -9\cdot 12a_0\cdot\eta_3-9\cdot 3a_1\cdot\eta_2+3\cdot 2a_2\cdot \eta_1,\\
		\gamma_2 &=& -9\cdot 6a_1\cdot\eta_3-9\cdot 4a_2\cdot\eta_2+3\cdot 6a_3\cdot\eta_1, \\
		\gamma_3 &=& -9\cdot 2a_2\cdot\eta_3-9\cdot 3a_3\cdot\eta_2+3\cdot 12a_4\cdot\eta_1.
	\end{eqnarray*}
	Restricting to monic polynomials by setting $a_0=1$ and making a few simplifications, we get a minimal presentation
	\[\Gamma(U_y,\Scal)=\coker\bgroup\begin{gmatrix}
		3(2a_3-a_1a_2) & 4a_4-a_2^2 \\
		\gamma_3-\frac{1}{2}a_1\gamma_2+\frac{3}{16}a_1^2\gamma_1 & \frac{1}{8}a_3\gamma_1-\frac{1}{6}a_2\gamma_2+\frac{1}{16}a_1a_2\gamma_1
	\end{gmatrix}\egroup.\]
\end{example}

\section{\label{SecArnold}The construction of Arnol'd}

The original definition of the open swallowtail is was given by Arnol'd in \cite{Arn81}. In this section, we describe his construction and show that it is equivalent to Definition \ref{DefOpenSwallowtail}.

To do so, we restrict to the affine subset $\{y\neq 0\ \mbox{and}\ a_0\neq 0\}$ of $\Pbf$. This may be identified with the space of monic polynomials of degree $n$ in one variable $x$. It will be convenient to work in different coordinates, namely, those associated to the \textit{divided powers} of $x$. For $1\leq i\leq n$, set $s_i:=(n!/(n-i)!)a_i$. Then $s_1,\dots,s_n$ are identified with the coefficients of the polynomial
\[ x^{(n)}+s_1x^{(n-1)}+s_2x^{(n-2)}+\cdots+s_n, \]
where $x^{(k)}:=x^k/k!$ is the $k$th \textit{divided power} of $x$. For $n\geq k>0$, differentiation of a polynomial in $x$ with respect to $x$ defines a finite map $\Sigma_{n,k}\to\Sigma_{n-1,k-1}$, which in turn defines a tower of varieties terminating at the discriminant $\Delta_n=\Sigma_{n,2}$. With respect to the coordinates $s_1,\dots,s_n$, this map is just projection $(s_1,\dots,s_n)\mapsto(s_1,\dots,s_{n-1})$.

For $n\geq k\geq 2$, we define the $k$th caustic $\Sigma_{n,k}$ to be the locus of polynomials of degree $n$ with a root of multiplicity at least $k$. Differentiation of such a polynomial with respect to the indeterminate gives rise to a polynomial of degree $n-1$ and a root of multiplicity $k$. This defines, for $n\geq k>2$, a finite map $\Sigma_{n,k}\to\Sigma_{n-1,k-1}$. Repeating this process defines a \textit{tower of caustics}
\begin{equation}\label{EqTowerOfCaustics}
	\cdots \to \Sigma_{n+i+1,i+3} \to \Sigma_{n+i,i+2} \to \Sigma_{n+i-1,i+1} \to \cdots \to \Sigma_{n,2} = \Delta_n
\end{equation}
terminating at the (affine) classical discriminant. Since each of the maps is finite and birational, the varieties in this tower share a common normalization.

Givental shows in \cite{Giv82} that this tower stabilizes, as in the following proposition.

\begin{proposition}[\cite{Giv82}, Theorem 1]
	For $n\geq 2$ and for $i>n-3$, the differentiation map $\Sigma_{n+i,i+2}\to\Sigma_{n+i-1,i+1}$ is an algebraic isomorphism. Thus the tower \eqref{EqTowerOfCaustics} stabilizes at $i=n-3$.
\end{proposition}

This motivates the following definition.

\begin{definition}
	The \textit{(geometric) open swallowtail} associated to $\Delta_n$, denoted $\Sigma_n$ (or $\Sigma$ when $n$ is understood), is the variety obtained at the point where the tower \eqref{EqTowerOfCaustics} stabilizes, namely, $\Sigma_{2n-3,n-1}$ in the notation above.
\end{definition}

\begin{figure}
	\begin{center}
		\includegraphics[width=6cm]{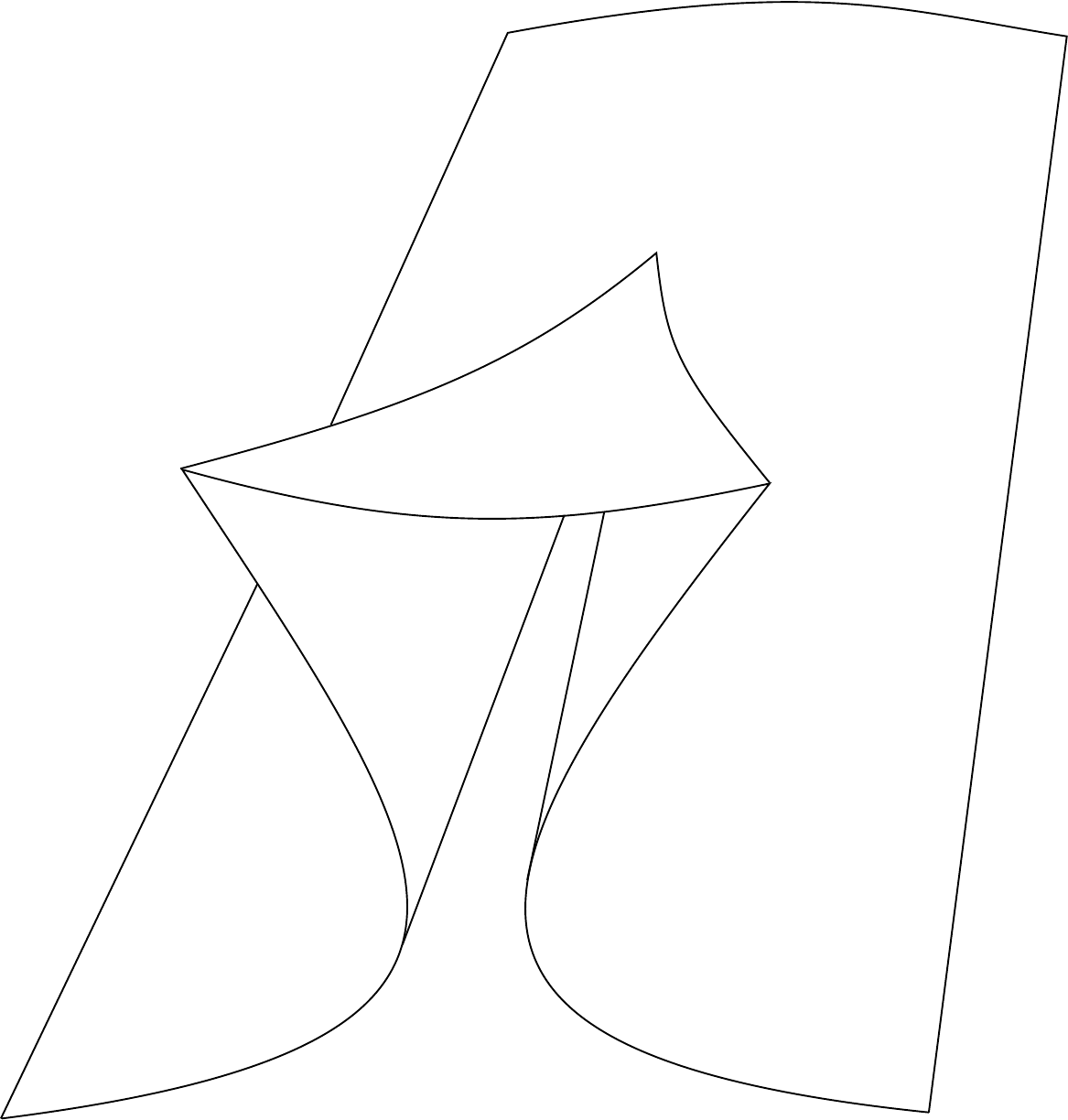}
	\end{center}
	\caption{\label{FigOpenSwallowtail}A schematic diagram of the open swallowtail $\Sigma_4$}
\end{figure}

The term \textit{open swallowtail} is used because $\Sigma$ is a partial normalization of $\Delta$ in which the self-intersection locus bifurcates, but the caustic remains. This is intuitively clear: a polynomial of degree $2n-3$ cannot have two roots of multiplicity $n-1$. However, it certainly can have a root of multiplicity strictly greater than $n-1$.

Since the differentiation map is finite, the coordinate ring $\Ocal_\Sigma$ of $\Sigma$ is a finite module over $\Ocal_\Delta$ which embeds in $\Ocal_{\bar{\Delta}}$. As such, it agrees with the algebraic open swallowtail introduced above. More precisely, the module $\Gamma(U_y,\Scal)$, to which we shall refer again as $\Scal$, is a module over the affine discriminant $\Ocal_\Delta$ which also by definition embeds in $\Ocal_{\bar{\Delta}}$. We have the following result.

\begin{theorem}
	\label{ThmArnoldsOS}
	The image of $\Ocal_\Sigma$ in $\Ocal_{\bar{\Delta}}$ equals $\Scal$.
\end{theorem}

To prove Theorem \ref{ThmArnoldsOS}, we need the following technical lemma of Givental which shows how the functions $s_{n-2+i}$ on $\Sigma_{n+k,k+2}$ for $k\geq i>0$ embed in the normalization $\bar{\Sigma}_{n+k,k+2}\cong\bar{\Delta}_n$.

\begin{lemma}[\cite{Giv82}, Lemma 2]
	\label{LemmaCharImageSigma}
	For $i>0$, $s_{n-2+i}=\pm\int_0^xf''(t)t^{(i-1)}\ dt$.
\end{lemma}

\begin{lemma}
	\label{PropDiffIsom}
	For $n\geq 2$ and $k\geq 0$, $\Omega^1_{\bar{\Sigma}_{n+k,k+2}/\Sigma_{n+k,k+2}}\cong\Omega^1_{\bar{\Sigma}_{n,2}/\Sigma_{n,2}}$.
\end{lemma}

\begin{proof}
	The modules $\Omega^1_{\bar{\Sigma}_{n+k,k+2}/\Sigma_{n+k,k+2}}$ and $\Omega^1_{\bar{\Sigma}_{n,2}/\Sigma_{n,2}}$ are the cokernels of the Jacobian matrices of the normalization maps $\pi:\bar{\Sigma}_{n+k,k+2}\to\Sigma_{n+k,k+2}$ and $\pi':\bar{\Sigma}_{n,2}\to\Sigma_{n,2}$. We use the local coordinates $x,s_1,\dots,s_{n-2}$ for $\bar{\Sigma}_{n+k,k+2}\cong\bar{\Sigma}_{n,2}$ and local coordinates $s_1,\dots,s_{n+k}$, respectively $s_1,\dots,s_n$, for $\Sigma_{n+k,k+2}$, respectively $\Sigma_{n,2}$. It suffices to show that, for $i>1$, the form $ds_{n+i-2}$ is a $\bar{\Sigma}_{n,2}$-linear combination of $ds_1,\dots,ds_{n-1}$. For $i>0$, differentiating the formula for $s_{n+i-2}$ given in Lemma \ref{LemmaCharImageSigma} with respect to the local coordinates on $\bar{\Sigma}_{n,2}$, we obtain the form
	\[ ds_{n+i-2} = \pm\left(f''(x)x^{(i-1)}\ dx + \binom{n+i-4}{i-1}x^{(n+i-3)}\ ds_1 + \cdots + \binom{i-1}{i-1}x^{(i)}\ ds_{n-2}\right). \]
	Setting $i=1$, we obtain the form
	\[ ds_{n-1} = \pm\left(f''(x)\ dx + x^{(n-2)}\ ds_1 + \cdots + x\ ds_{n-2}\right). \]
	Therefore, for $i>1$, $ds_{n+i-2}$ is the sum of $\pm x^{(i-1)}ds_{n-1}$ and a suitable $\bar{\Sigma}_{n,2}$-linear combination of $ds_1,\dots,ds_{n-2}$. The claim follows.
\end{proof}

\begin{proof}[Proof of Theorem \ref{ThmArnoldsOS}]
	Lemma \ref{PropDiffIsom} shows that $\Omega^1_{\bar{\Sigma}/\Sigma}\cong\Omega^1_{\bar{\Delta}/\Delta}$, which implies that $\Ocal_\Sigma$ is contained in $\Scal$. We claim that the embedding $\Ocal_\Sigma\hookrightarrow\Scal$ is surjective. From Theorem \ref{ThmCharOmega1}, we see that the kernel of $d$ consists of those elements $g\in\Ocal_{\bar{\Delta}}$ which, after being differentiated with respect to $x$, are divisible by $f''$. Such elements are of the form
	\[g(x,s_1,\dots,s_{n-2})=\int_0^xh_1(t,s_1,\dots,s_{n-2})f''\ dt+h_2(s_1,\dots,s_{n-2}).\]
	Writing $h_1$ as a polynomial in $t$ with coefficients in $\KK[s_1,\dots,s_{n-2}]$, we see that $g-h_2$ is an $\Ocal_\Delta$-linear combination of elements of the form
	\[ \int_0^xf''(t)t^{(i-1)}\ dt, \quad i\geq 1 \]
	and therefore, by Lemma \ref{LemmaCharImageSigma}, in the image of $\Ocal_\Sigma$. On the other hand, $h_2$ is in the image of $\Ocal_\Delta$ and a fortiori of $\Ocal_\Sigma$. This proves the claim.
\end{proof}

In \cite{Sev04} Sevenheck and van Straten point out that Givental's results in \cite{Giv88} imply that $\Sigma$ is a Cohen-Macaulay variety, and hence that $\Ocal_\Sigma$ is maximal Cohen-Macaulay over $\Ocal_\Delta$. Our results yield a new, algebraic proof that the open swallowtail of Arnol'd is Cohen-Macaulay.

\begin{corollary}
	The open swallowtail $\Sigma$ is Cohen-Macaulay.
\end{corollary}

\begin{proof}
	Combine Theorem \ref{ThmArnoldsOS} with Proposition \ref{PropOSProps}.
\end{proof}

The following result gives a generating set for $\Ocal_\Sigma$ as an $\Ocal_\Delta$-module. It was first proved by Givental in \cite{Giv88}, but we give an alternative proof using the tools developed here.

\begin{theorem}
	\label{ThmOSGens}
	The module $\Ocal_\Sigma$ is generated over $\Ocal_\Delta$ by 1 and $s_{n+1},\dots,s_{2n-3}$.
\end{theorem}

\begin{proof}
	Clearly $\Ocal_\Sigma$ is generated as an $\Ocal_\Delta$-algebra by the aforementioned elements. From Theorem \ref{ThmArnoldsOS} and Proposition \ref{PropMinGenerators}, $\Ocal_\Sigma$ as an $\Ocal_\Delta$-module has $n-2$ minimal generators. In order to prove the claim, it suffices to show that no element $s_{n+i}$ lies in the submodule generated by the remaining generators $1,s_{n+1},\dots,\hat{s}_{n+i},\dots,s_{2n-3}$. Suppose by way of contradiction otherwise, that
	\[ s_{n+i}=\sum_{j=1,\dots,i-1,i+1,\dots,n}f_js_{n+j}. \]
	Since $\Ocal_\Sigma$ is graded, we may take the above equation to be quasihomogeneous. Since each $s_{n+i}$ has degree $n+i$, degree considerations indicate that $f_j=0$ for $j>i$, that is
	\[ s_{n+i}=\sum_{j=1,\dots,i-1}f_js_{n+j}. \]
	But that would imply that the map $\Sigma_{n+i,i+2}\to\Sigma_{n+i-1,i+1}$ were an isomorphism. However, since, for $i\leq n-3$, each map $\Sigma_{n+i,i+2}\to\Sigma_{n+i-1,i+1}$ reduces the dimension of the self-intersection locus, these maps are not isomorphisms. This is a contradiction, and the claim follows.
\end{proof}

\section{\label{SecConductor}The conductor of the open swallowtail}

We now show that the map $\Sigma\to\Delta$ is an isomorphism outside of self-intersection locus of $\Delta$. More precisely, we prove the following result.

\begin{theorem}
	\label{ThmConductorOS}
	Let $\cfrak$ be the conductor of the map $\Sigma\to\Delta$. Then the closed subset of $\Delta$ defined by $\cfrak$ is the self-intersection locus of $\Delta$.
\end{theorem}

We start with the following immediate corollary of Theorem \ref{ThmSingSigma}.

\begin{proposition}
	\label{PropJacODBar}
	The closed subset of $\Delta$ defined by the conductor $\dfrak$ of the normalization map $\Ocal_\Delta\hookrightarrow\Ocal_{\bar{\Delta}}$ is the singular locus of $\Ocal_\Delta$. It has two irreducible components: the caustic and the self-intersection locus.
\end{proposition}

We require the following classical result, which describes the behaviour of the discriminant at sufficiently generic points of the caustic.

\begin{lemma}
	\label{ThmGenPtCaustic}
	Let $p\in\Gamma$ be a point which is not on the self-intersection locus of $\Delta_n$. Then the germ of $\Delta_n$ at $p$ is formally isomorphic to the product of the cuspidal cubic $\Spec\KK[x,y]/(x^3-y^2)$ and a smooth factor.
\end{lemma}

\begin{lemma}
	\label{LemmaJacGenRed}
	The conductor $\dfrak$ is reduced at points of the caustic which are not in the self-intersection locus.
\end{lemma}

\begin{proof}
	Such points correspond to polynomials of degree $n$ with exactly one root of multiplicity exactly three and with all other roots distinct. Lemma \ref{ThmGenPtCaustic} implies that, locally at such a point, $\Delta$ is the product of a smooth factor and the cuspidal cubic. The claim now follows from the equivalent claim for the cuspidal cubic, which is an easy calculation.
\end{proof}

The following lemma indicates that, in the affine setting, $\Omega^1_{\bar{\Delta}/\Delta}$ is a Gorenstein module over $\Ocal_\Delta$.

\begin{lemma}
	\label{LemmaOmegaDuality}
	We have $\Ext^1_{\Ocal_\Delta}(\Omega^1_{\bar{\Delta}/\Delta},\Ocal_\Delta)\cong\Omega^1_{\bar{\Delta}/\Delta}$.
\end{lemma}

\begin{proof}
	Proposition \ref{PropGammaBarCM} and Theorem \ref{ThmCharOmega1} imply that $\Omega^1_{\bar{\Delta}/\Delta}$ is presented as a quotient of $\Ocal_{\bar{\Delta}}$ by a principal ideal thereof. That is, we have an exact sequence
	\[ \begin{diagram} 0 & \rTo & \Ocal_{\bar{\Delta}} & \rTo & \Ocal_{\bar{\Delta}} & \rTo & \Omega^1_{\bar{\Delta}/\Delta} & \rTo & 0 \end{diagram} \]
	of $\Ocal_{\bar{\Delta}}$-modules. Treating it as an exact sequence of $\Ocal_\Delta$ modules and applying $\Hom_{\Ocal_\Delta}(-,\Ocal_\Delta)$, we obtain the sequence
	\[ \begin{diagram} 0 & \rTo & \Hom_{\Ocal_\Delta}(\Ocal_{\bar{\Delta}},\Ocal_\Delta) & \rTo_\eta & \Hom_{\Ocal_\Delta}(\Ocal_{\bar{\Delta}},\Ocal_\Delta) & \rTo & \Ext^1_{\Ocal_\Delta}\left(\Omega^1_{\bar{\Delta}/\Delta},\Ocal_\Delta\right) \\ & \rTo & \Ext^1_{\Ocal_\Delta}(\Ocal_{\bar{\Delta}},\Ocal_\Delta) = 0, \end{diagram} \]
	where the vanishing on the right follows from $\Ocal_{\bar{\Delta}}$ being maximal Cohen-Macaulay over the Gorenstein ring $\Ocal_\Delta$. Since $\Ocal_{\bar{\Delta}}$ is a maximal Cohen-Macaulay module on a hypersurface and is presented by a symmetric matrix, it is self-dual. The map $\eta$ is just multiplication by the same non-zerodivisor which presents $\Omega^1_{\bar{\Delta}/\Delta}$. Thus both $\Omega^1_{\bar{\Delta}/\Delta}$ and $\Ext^1_{\Ocal_\Delta}(\Omega^1_{\bar{\Delta}/\Delta},\Ocal_\Delta)$ have the same presentations as $\Ocal_\Delta$-modules, and the claim follows.
\end{proof}

We shall require the following technical lemma, whose easy proof is left to the reader.

\begin{lemma}
	\label{LemmaCharConductor}
	Suppose that $X\to Y$ is a finite, birational map of irreducible affine varieties with conductor $\cfrak$. Then the map $\Hom_{\Ocal_Y}(\Ocal_X,\Ocal_Y)\to\Ocal_Y$ sending $\alpha$ to $\alpha(1)$ is an isomorphism of $\Hom_{\Ocal_Y}(\Ocal_X,\Ocal_Y)$ onto $\cfrak$.
\end{lemma}

\begin{proof}[Proof of Theorem \ref{ThmConductorOS}]
	We begin with the short exact sequence
	\[ \begin{diagram} 0 & \rTo & \Ocal_\Sigma & \rTo & \Ocal_{\bar{\Delta}} & \rTo^d & \Omega^1_{\bar{\Delta}/\Delta} & \rTo & 0 \end{diagram} \]
	which results from combining the definition of the open swallowtail $\Scal$ with Theorem \ref{ThmArnoldsOS}.	Applying $\Hom_{\Ocal_\Delta}(-,\Ocal_\Delta)$, we obtain
	\begin{equation}
		\label{EqHomExt}
		\begin{diagram}
			0 & \rTo & \Hom_{\Ocal_\Delta}(\Ocal_{\bar{\Delta}},\Ocal_\Delta) & \rTo^j & \Hom_{\Ocal_\Delta}(\Ocal_\Sigma,\Ocal_\Delta) & \rTo & \Ext^1_{\Ocal_\Delta}\left(\Omega^1_{\bar{\Delta}/\Delta},\Ocal_\Delta\right) & \rTo & 0
		\end{diagram},
	\end{equation}
	where exactness on the right follows from $\Ocal_{\bar{\Delta}}$ being maximal Cohen-Macaulay and $\Ocal_\Delta$ being Gorenstein. Now let $\cfrak$ be the conductor of $\Ocal_\Delta\hookrightarrow\Ocal_\Sigma$ and $\dfrak$ be the conductor of $\Ocal_\Delta\hookrightarrow\Ocal_{\bar{\Delta}}$. Lemma \ref{LemmaCharConductor} implies that $\cfrak\cong\Hom_{\Ocal_\Delta}(\Ocal_\Sigma,\Ocal_\Delta)$ and $\dfrak\cong\Hom_{\Ocal_\Delta}(\Ocal_{\bar{\Delta}},\Ocal_\Delta)$ and that the following diagram commutes:
	\[ \begin{diagram}
		\dfrak & \rInto^i & \cfrak \\
		\dEquals_\sim & & \dEquals_\sim \\
		\Hom_{\Ocal_\Delta}(\Ocal_{\bar{\Delta}},\Ocal_\Delta) & \rTo^j & \Hom_{\Ocal_\Delta}(\Ocal_\Sigma,\Ocal_\Delta),
	\end{diagram} \]
	where $i:\dfrak\to\cfrak$ is the natural inclusion and $j$ is the same map as given in \eqref{EqHomExt}.
	Applying these identifications, Lemma \ref{LemmaOmegaDuality}, and Theorem \ref{ThmCharOmega1} to \eqref{EqHomExt}, we obtain a commutative diagram
	\begin{equation}
		\label{EqCondDiagram}
		\begin{diagram}[height=.5cm]
			& & 0 & & 0 \\
			& & \dTo & & \dTo \\
			0 & \rTo & \dfrak & \rTo^i & \cfrak & \rTo & \Ocal_{\bar{\Gamma}} & \rTo & 0 \\ \\
			& & \dTo & & \dTo \\ \\
			& & \Ocal_\Delta & \rEquals & \Ocal_\Delta \\ \\
			& & \dTo & & \dTo \\ \\
			& & \Ocal_\Delta/\dfrak & \rTo & \Ocal_\Delta/\cfrak \\ \\
			& & \dTo & & \dTo \\
			& & 0 & & 0.
		\end{diagram}
	\end{equation}

	The map $i$ in \eqref{EqCondDiagram} is an isomorphism outside of the caustic $\Gamma$. Thus, outside of $\Gamma$, the support of $\cfrak$ equals that of $\dfrak$. Proposition \ref{PropJacODBar} implies that this support is precisely the self-intersection locus.

	It remains to show that, among points on $\Gamma$, $\cfrak$ is supported only at the intersection of $\Gamma$ and the self-intersection locus. The snake lemma applied to \eqref{EqCondDiagram} implies the existence of an exact sequence
	\[ \begin{diagram}
		0 & \rTo & \Ocal_{\bar{\Gamma}} & \rTo & \Ocal_\Delta/\dfrak & \rTo & \Ocal_\Delta/\cfrak & \rTo & 0.
	\end{diagram} \]
	Let $p$ be a point of $\Gamma$ not in the self-intersection locus. Then, by Theorem \ref{ThmSingSigma}, $p$ is a smooth point of $\Gamma$, so the map $\bar{\Gamma}\to\Gamma$ is an isomorphism at $p$. Also, at $p$, by Lemma \ref{LemmaJacGenRed}, $\Ocal_\Delta/\dfrak\cong\Ocal_\Gamma$. Thus $\Ocal_\Delta/\cfrak$ is not supported at $p$. This proves the claim.
\end{proof}

\section{\label{SecApplication}Application to the root structure of a univariate polynomial}

We now give a compelling application of the matrix of the open swallowtail. The rank of a minimal presentation $A$ of the normalization of $\Delta_n$, specialized to a particular polynomial $f$, is the number of distinct roots of $f$ minus one. This follows from the description of the normalization of $\Delta$ given in Section \ref{SecPreliminaries}. The nullity of the specialization of $A$ to $f$ is the length of the fibre of the normalization map over $f$. The fibre of $f$ consists of pairs $(f,x)$ where $x$ is a repeated root of $f$. The multiplicity of the point $(f,x)$ in the fibre is, in turn, one less than the multiplicity of $x$ as a root of $f$.

However, the rank of the matrix of the normalization is not able to distinguish between degenerate repeated roots and multiple repeated roots. For example, it cannot detect whether a polynomial has $n-2$ distinct roots because it has two distinct pairs of double roots or because it has one root of multiplicity 3. We show here that the matrix of the open swallowtail can make this distinction.

\begin{theorem}
	\label{ThmOSRanks}
	Let $B$ be a minimal presentation of $\Ocal_\Sigma$ over the ambient ring of $\Ocal_\Delta$. Suppose $B$ is specialized to some polynomial $f(x)=x^n+a_1x^{n-1}+\cdots+a_n$. The nullity of the resulting matrix is at least 2 if and only if $f(x)$ has more than one distinct pair of double roots.
\end{theorem}

\begin{proof}
	Fix $i\geq 0$. The ideal defining the locus of polynomials $f(x)$ for which the specialized matrix $B$ has nullity at least $2$ is the radical of the Fitting ideal $F_1(\Sigma)$. By \cite{Eis95}, Proposition 20.6, the points defined by $F_1(\Sigma)$ are those where $\Ocal_\Sigma$ cannot be generated by one element, that is, where the map $\Sigma\to\Delta$ is not an isomorphism. The locus of such points is precisely the zero locus of the conductor of the map $\Sigma\to\Delta$, which, by Theorem \ref{ThmConductorOS}, is precisely the self-intersection locus, whence the claim.
\end{proof}

It is natural to ask what would be the meaning of the nullity of the matrix of the open swallowtail being strictly greater than two. It is likely that, for $i\geq 0$, the nullity of this matrix is at least $i+2$ if and only if the polynomial $f$ is the $i$th derivative of a polynomial with more than one root of multiplicity $i+2$. However, we do not have a proof of this.

\bibliographystyle{amsalpha}
\bibliography{open-swallowtail}

\end{document}